\documentclass[twocolumn,amssymb,aps,nofootinbib,secnumarabic]{revtex4}
\usepackage{times,mathptmx,amsthm2000,amsmath2000}
\usepackage{pstricks,pst-node}
\newtheorem{theorem}{Theorem}
\newtheorem{proposition}[theorem]{Proposition}
\newtheorem{corollary}[theorem]{Corollary}
\theoremstyle{remark}
\newtheorem*{remark}{Remark}

\newcommand{\Z}{\mathbb{Z}}
\newcommand{\R}{\mathbb{R}}

\newcommand{\braket}[1]{{\langle #1 \rangle}}
\renewcommand{\d}{\partial}

\newcommand{\cA}{\mathcal{A}}
\newcommand{\cB}{\mathcal{B}}
\newcommand{\cC}{\mathcal{C}}
\newcommand{\cD}{\mathcal{D}}
\newcommand{\bcD}{\mathcal{\overline{D}}}
\newcommand{\cI}{\mathcal{I}}
\newcommand{\cT}{\mathcal{T}}

\newcommand{\eatline}{\vspace{-\baselineskip}}
\newcommand{\eq}[2]{\begin{equation}\label{#1}#2\end{equation}}

\newenvironment{fullfigure}[2]
    {\begin{figure}[htb]\def\ffa{#1}\def\ffb{#2}}
    {\vspace{\baselineskip}\caption{\ffb.}\label{\ffa}\end{figure}}

\newcommand{\fig}[1]{Figure~\ref{#1}}
\newcommand{\thm}[1]{Theorem~\ref{#1}}
\newcommand{\cor}[1]{Corollary~\ref{#1}}
\renewcommand{\sec}[1]{Section~\ref{#1}}

\psset{linewidth=.5pt,dash=4pt 4pt,unit=.25in,arrowsize=2pt 3}
\SpecialCoor

\begin{document}
\title{A generalization of Filliman duality}
\author{Greg Kuperberg}
\thanks{Supported by NSF grant DMS \#0072342}
\affiliation{UC Davis}
\email{greg@math.ucdavis.edu}
\begin{abstract}
Filliman duality expresses (the characteristic measure of) a convex polytope
$P$ containing the origin as an alternating sum of simplices that share
supporting hyperplanes with $P$.  The terms in the alternating sum are given by
a triangulation of the polar body $P^\circ$.  The duality can lead to useful
formulas for the volume of $P$.  A limiting case called Lawrence's algorithm
can be used to compute the Fourier transform of $P$.

In this note we extend Filliman duality to an involution on the space of
polytopal measures on a finite-dimensional vector space, excluding polytopes
that have a supporting hyperplane coplanar with the origin.  As a special case,
if $P$ is a convex polytope containing the origin, any realization of $P^\circ$
as a linear combination of simplices leads to a dual realization of $P$.
\end{abstract}
\maketitle

\section{Introduction}

If $P \subset \R^d$ is a polytopal region, let $[P]$ denote
the restriction of Lebesgue measure to $P$.  The measure $[P]$ can
also be called the characteristic measure, by analogy with the characteristic
function.  We consider measures rather than functions so that if $P$
and $Q$ are two regions with disjoint interiors, then
$$[P \cup Q] = [P] + [Q]$$
even if $P$ and $Q$ are not disjoint at the boundary.  If $P$ is convex and the
origin lies in its interior, then $P$ admits a polar body $P^\circ$. Say that a
polytope $P$, convex or not, is \emph{codegenerate} if one of its facets is
coplanar with the origin.  Every non-codegenerate simplex $\Delta$ admits a
polar simplex $\Delta^\circ$. Let $\cT$ be a triangulation of $P^\circ$ by
non-codegenerate simplices. In this circumstance Filliman \cite{Filliman:duals}
showed that
$$[P] = \sum_{\Delta \in \cT} (-1)^{\sigma(\Delta)} [\Delta^\circ],$$
where $\sigma(\Delta)$ is a certain sign function.  This formula
is called \emph{Filliman duality}.  \fig{f:pentagon} shows
an example:  a triangulation of a pentagon $P$ and the dual realization of
$P^\circ$ as a triangle with two smaller triangles subtracted.

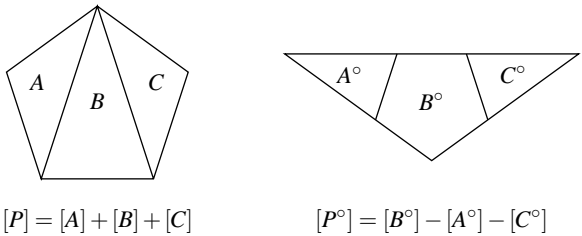
\begin{fullfigure}{f:pentagon}{A triangulation of a pentagon
    and its Filliman dual}
\pspicture(-2,-3)(2,2)
\psline(2;234)(2;306)(2;18)(2;90)(2;162)(2;234)(2;90)(2;306)
\rput(1.309;162){$A$}
\rput(0,0){$B$}
\rput(1.309;18){$C$}
\rput(0,-2.5){$[P] = [A] + [B] + [C]$}
\endpspicture
\pspicture(-5,-3)(3,2)
\pspolygon(3.236;18)(3.236;162)(1.236;270)
\psline(1.236;198)(1.236;126)
\psline(1.236;342)(1.236;54)
\rput(1.8;162){$A^\circ$}
\rput(0,0){$B^\circ$}
\rput(1.8;18){$C^\circ$}
\rput(0,-2.5){$[P^\circ] = [B^\circ] - [A^\circ] - [C^\circ]$}
\endpspicture
\end{fullfigure}

Two special cases of Filliman duality are notable.  First, if each simplex
$\Delta \in \cT$ shares vertices with $P^\circ$, then $\Delta^\circ$ shares
supporting hyperplanes with $P$.  If $P$ has few vertices and many sides, then
it is reasonable to compute its volume as the sum of the volumes of the
simplices in a triangulation.  But if $P$ has many vertices and few sides, it
is more efficient to use a triangulation of $P^\circ$ via Filliman duality. 
Second, if $\cT$ is the cone of  a triangulation of $\d P^\circ$, then as the
apex of the cone converges to the origin, each dual simplex $\Delta^\circ$ with
$\Delta \in \cT$ converges to an affine orthant emanating from a vertex of $P$.
Thus we can express $[P]$ as an alternating sum of such orthants.  This
limiting case of Filliman duality is called Lawrence's algorithm
\cite{Lawrence:computation}. It is useful not only for finding the volume of
$P$, but also for computing the Fourier transform of $[P]$.

In this note we extend Filliman duality to an involution on non-codegenerate,
integral, polytope measures on $\R^d$.  More precisely, let $\cA$ be the
abelian group of signed measures $\mu$ on $\R^d$ of the form
$$\mu = \sum_{i=1}^n \alpha_i [P_i],$$
where each $P_i$ is a non-codegenerate polytope and $\alpha_i \in \Z$. If
$\Delta$ is a non-codegenerate simplex, let $\sigma(\Delta)$ be the number
of supporting hyperplanes of $\Delta$ that separate it from the
origin.  Let $\cD \subset \cA$ consist of measures of the form
$[\Delta]$.

\begin{theorem} The involution $\Phi:\cD \to \cD$
defined by
$$\Phi([\Delta]) = (-1)^{\sigma(\Delta)} [\Delta^\circ]$$
extends uniquely to an automorphism $\Phi:\cA \to \cA$.
\label{th:main}\end{theorem}

The following corollary captures the original Filliman duality
as part of the involution $\Phi$:

\begin{corollary}
If $P$ is a convex polytope containing the origin in its interior
and
$$[P] = \sum_{i=1}^n \alpha_i [\Delta_i]$$
for simplices $\Delta_1,\dots,\Delta_n$, then
$$[P^\circ] = \sum_{i=1}^n (-1)^{\sigma(\Delta)}\alpha_i [\Delta_i^\circ].$$
\label{c:convex}\end{corollary}

Note that our involution $\Phi$ is not the same as the Euler involution on the
polytope algebra defined by McMullen \cite{McMullen:polytope}. Nonetheless, as
the group $\cA$ and the map $\Phi$ on it are ultimately a disguised
specialization of a known involution in valuation theory, the polarity map on
cones (see \sec{s:extra}).  Thus the real significance of \thm{th:main} and its
proof is not that it describes an essentially new involution, as the author
once thought, but rather that it relates three distinct constructions in
combinatorial geometry: Filliman duality, valuation theory, and the stellar
subdivision theorem.

\section{Proof of \thm{th:main}}

In this section we assume that all simplices and other polytopes are
non-codegenerate except where noted.

We will consider signed polytopes in order to absorb the sign that appears in
duality for simplices.  A signed polytope $P$ is a polytope together with a
formal sign, either $+$ or $-$.  Characteristic measures on signed polytopes
are defined by the rule  $[-P] = -[P]$.  If $\Delta \subset \R^d$ is a
simplex with positive sign, we define
$$\Delta^* = (-1)^{\sigma(\Delta)}\Delta^\circ
\qquad (-\Delta)^* = -(\Delta^*).$$
Also we recall the definition of $\Delta^\circ$.  If $\Delta$ has
vertices $v_0,\dots,v_d$, then $\Delta^\circ$ is bounded by
the hyperplanes $H_{v_0},\dots,H_{v_d}$, where for any vector $v$,
$H_v$ is defined as
$$H_v = \left\{ w | \braket{w,v} = 1 \right\}.$$

Let $\cB$ be the abelian group freely generated by simplex measures
$[\Delta]$.  By summing the terms of each element in $\cB$, we obtain a
homomorphism $\pi:\cB \to \cA$.  It is surjective because every polytopal
region can be tiled by simplices.  The involution $\Phi$ extends tautologically
to $\cB$.  \thm{th:main} then asserts that $\Phi$ preserves $\ker \pi$. In
order to prove this we first give a characterization of the kernel.
The characterization depends on a refinement of the stellar subdivision
due to M. H. A. Newman \cite{Lickorish:complexes,Newman:theorem}.

\begin{theorem}[Newman] Any two triangulations of a polytope $P \subset \R^d$
are equivalent under stellar moves applied on edges.
\label{th:stellar}
\end{theorem}

\begin{corollary} The kernel $\ker \pi$ is generated by the
relators
$$[\Delta] - [\Delta_1] - [\Delta_2],$$
where $\Delta$ is a simplex tiled by two simplices 
$\Delta_1$ and $\Delta_2$.
\label{c:split}
\end{corollary}

\begin{proof} Let $\cI$ be the subgroup of $\cB$ generated by the relators. 
Clearly $\cI \subseteq \ker \pi$; we wish to show that $\ker \pi \subseteq
\cI$. Assume a general linear dependence of simplices
\eq{e:dep}{\sum_{i=1}^n \alpha_i [\Delta_i] = 0}
in $\cA$.  Equivalently, in $\cB$,
$$\sum_{i=1}^n \alpha_i [\Delta_i] \in \ker \pi.$$
The union of the simplices,
$$P = \cup_{i=1}^n \Delta_i,$$
is a compact polytopal region in $\R^d$.  It admits a triangulation $\cT$ that
refines each simplex $\Delta_i$.  Let $\cT_i$ be the restriction of $\cT$ to
the simplex $\Delta_i$.  By \thm{th:stellar}, the triangulation $\cT_i$ can be
obtained from the tautological triangulation of $\Delta_i$ by itself by stellar
moves applied to edges.  A stellar move on some edges $e$ can be effected by
dividing each simplex containing $e$ into two simplices, the geometric move
captured by the relator.  Therefore 
$$[\Delta_i] - \sum_{\Delta \in \cT_i} [\Delta] \in \cI$$
in $\cB$. At the same time, equation \eqref{e:dep} implies that
$$\sum_{i=1}^n \alpha_i \sum_{\Delta \in \cT_i} [\Delta] = 0$$
in $\cB$, since each simplex in $\cT$ must be covered a total of 0 times.
Therefore $\ker \pi \subseteq \cI$, as desired.
\end{proof}

\begin{remark} Call the move of dividing a simplex into two an \emph{elementary
dissection}.  An interesting fact closely related to \cor{c:split} is that any
two simplicial dissections of a polytopal region are connected by elementary
dissections.
\end{remark}

\begin{proof}[Proof of \thm{th:main}]
In light of \cor{c:split}, we only need to check that $\Phi$ preserves an
elementary dissection of a simplex.  If $v_0,v_1,\dots,v_d$ are 
the vertices of a simplex $\Delta$, we let the word
$v_0v_1\dots v_d$ denote $\Delta$ if
$$v_1-v_0,v_2-v_0,\dots,v_d-v_0$$
is a positive basis of $\R^d$, and otherwise we let it denote $-\Delta$.
If $w_i$ is the vertex of $\Delta^*$ opposite to the hyperplane
$H_{v_i}$, it follows that
$$\Delta^* = (-1)^d w_0w_1\dots w_d.$$

Suppose that $x_1$,$x_2$, and $x_3$ are 3
collinear points in $V$, and suppose that $v_0,\ldots,v_{d-2}$ are $d-1$ other
points affinely independent from any two of $x_0$, $x_1$, and $x_2$.
Then
\begin{multline}
[v_0\ldots v_{d-2}x_1x_2] + [v_0\ldots v_{d-2}x_2x_3] \\
    \hfill + [v_0\ldots v_{d-2}x_3x_1] = 0
\label{e:split1}\end{multline}
expresses an elementary dissection.  Applying $(-1)^d\Phi$ to both sides
produces
\begin{multline}
[w_0w_1\ldots w_{d-2}y_2y_1] + [w_0w_1\ldots w_{d-2}y_3y_2] \\
    \hfill + [w_0w_1\ldots w_{d-2}y_1y_3] = 0.
\label{e:split2}\end{multline}
Here each point $w_i$ lies in the hyperplanes $H_{v_j}$ for $i \ne j$ and 
in the hyperplanes $H_{x_j}$ for all $j$.  Each point $y_i$
lies in the hyperplane $H_{v_j}$ for all $j$ and in the hyperplane
$H_{x_i}$.  Evidently equation~\eqref{e:split2} is the same equation
up to sign as equation~\eqref{e:split1}.
\fig{f:split} shows an example.
\end{proof}

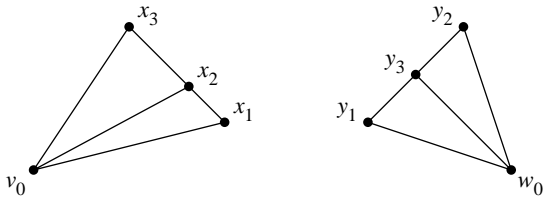
\begin{fullfigure}{f:split}{An elementary dissection and its dual}
\pspicture(-1,-1)(11,4)
\qdisk(0,0){.1} \qdisk(2,3){.1} \qdisk(3.25,1.75){.1} \qdisk(4,1){.1}
\psline(0,0)(2,3)(4,1)(0,0)(3.25,1.75)
\uput{.2}[30](2,3){$x_3$}
\uput{.2}[30](3.25,1.75){$x_2$}
\uput{.2}[30](4,1){$x_1$}
\uput{.2}[225](0,0){$v_0$}
\qdisk(10,0){.1} \qdisk(9,3){.1} \qdisk(8,2){.1} \qdisk(7,1){.1}
\psline(10,0)(9,3)(7,1)(10,0)(8,2)
\uput{.2}[315](10,0){$w_0$} \uput{.2}[150](9,3){$y_2$}
\uput{.2}[150](8,2){$y_3$} \uput{.2}[150](7,1){$y_1$}
\endpspicture
\eatline\end{fullfigure}

\section{Proof of \cor{c:convex} and examples}

\cor{c:convex} follows immediately from \thm{th:main} and the
following proposition.

\begin{proposition} If $P \subset \R^d$ is a convex polytope that strictly
contains the origin, then $\Phi([P]) = [P^*]$, the polar body of $P$.
\label{p:polar}
\end{proposition}

\begin{proof} We first assume that any $d$ vertices of $P$ are linearly
independent, or equivalently affinely independent from the origin.  Let $\cT$
be a triangulation of $P$ with no vertices in the interior of $P$.  We claim,
first, that any point $w$ in the interior of $P^*$ is covered by the dual of
exactly one simplex. A unique simplex $\Delta_0 \in \cT$ contains the origin.
The inclusions $0 \in \Delta_0 \subset P$ imply that $\Delta_0^* \supset P^*$,
so $\Delta^*$ covers $w$.  If $\Delta \in \cT$ is another simplex, then there
exists a vertex $v$ of $\Delta$ which is separated from the origin by the
opposite face of $\Delta$. It follows that $\Delta^*$ is separated from the
origin by the hyperplane $H_v$.  Since $v$ is also a vertex of $P$, $H_v$ is a
supporting hyperplane of $P^*$.  Therefore $H_v$ separates $\Delta^*$ from
$P^*$ and $\Delta^*$ does not contain $w$.  This establishes the first claim.

We claim, second, that if $w$ is in the exterior of $P^*$, there exists a
triangulation $T$ of $P$ such that $w$ is not covered by $\Delta^*$ for any
$\Delta \in T$.  There exists a vertex $v$ of $P$ such that the hyperplane
$H_v$ separates $w$ from $P^*$.  Let $\cT$ be a fan triangulation all of whose
simplices contain $v$.  If $\Delta \in \cT$, then $H_v$ separates $\Delta$ from
$w$, as desired.  The two claims together with \thm{th:main} establish the
proposition under the independence assumption on $P$.

Finally we assume that $P$ is arbitrary. The argument so far establishes the
proposition for the polytope $P-v$ for a dense set of vectors $v$ in the
interior of $P$.  Namely $v$ can be any point that does not lie on a hyperplane
affinely spanned by vertices of $P$.  But if $v$ is in the interior of $P$,
then $P-v$ has a non-codegenerate triangulation whether or not $v$ lies on such
a hyperplane.  It follows that $\Phi([P-v])$ varies continuously in a
neighborhood of $v$.  Thus the truth of the proposition for a dense set of $v$
implies its truth for all $v$ in the interior of $P$.  In particular the
proposition holds for $v=0$.
\end{proof}

We conclude with two examples of polygonal regions in the plane
and their images under $\Phi$.
\fig{f:square1} shows the dual of a square with corners
$(1,-1)$, $(3,-1)$, $(1,1)$, and $(3,1)$.  The dual is a negative
measure in the interior of a kite shape.  Finally 
\fig{f:square2} shows the dual of a square with corners
$(1,1)$, $(1,2)$, $(2,1)$, and $(2,2)$.  The dual is the difference
(in the sense of signed measures) between two triangles with disjoint
interiors.

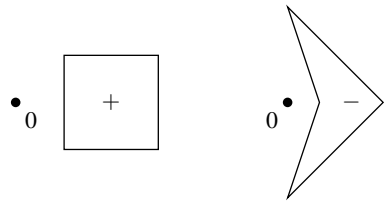
\begin{fullfigure}{f:square1}{A square offset to the right and its dual}
\pspicture(-.5,-2.5)(3.5,2.5)
\qdisk(0,0){.1}\rput[tl](.2,-.2){$0$}
\psframe(1,-1)(3,1)
\rput(2,0){$\mathbf{+}$}
\endpspicture
\hspace{1cm}
\pspicture(-.5,-2.5)(2.5,2.5)
\qdisk(0,0){.1}\rput[tr](-.2,-.2){$0$}
\pspolygon(.667,0)(0,2)(2,0)(0,-2)
\rput(1.333,0){$\mathbf{-}$}
\endpspicture
\end{fullfigure}

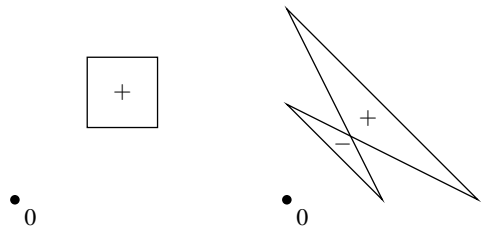
\begin{fullfigure}{f:square2}{A square offset diagonally and its dual}
\pspicture(-.5,-.5)(3.5,3.5)
\qdisk(0,0){.1}\rput[tl](.2,-.2){$0$}
\psframe(1.5,1.5)(3,3)
\rput(2.25,2.25){$\mathbf{+}$}
\endpspicture
\hspace{1cm}
\pspicture(-.5,-.5)(4.5,4.5)
\qdisk(0,0){.1}\rput[tl](.2,-.2){$0$}
\pspolygon(2,0)(0,4)(4,0)(0,2)
\rput(1.167,1.167){$\mathbf{-}$}
\rput(1.7,1.7){$\mathbf{+}$}
\endpspicture
\end{fullfigure}

\section{Involutions on cones}
\label{s:extra}

In response to the first version of this article, Alexander Barvinok suggested
that \thm{th:main} is related to the fact that the polarity involution $\Pi$ on
spherical convex polytopes extends to a valuation.  In this section we show, in
outline, that $\Phi$ becomes a restriction of $\Pi$ after suitably mapping
measures on $\R^d$ to functions on the sphere $S^d$.

Let $\cC$ be the abelian group of integer-valued functions on $S^d$ spanned by
the characteristic functions of closed (or equivalently open) spherical
polytopes of any dimension $\le d$.  If $P \subset S^d$ is closed and convex,
then it admits a polar dual
$$P^\circ = \{x | \forall y \in P, \braket{x,y} \ge 0 \},$$
where the inner product uses the defining embedding $S^d \subset \R^{d+1}$.
Also, if $P \subset S^d$, let $\chi_P$ be the characteristic function of $P$. 
The following result in combinatorial geometry is known but unattributed (see
Barvinok \cite{Barvinok:gsm} and Lawrence \cite{Lawrence:polarity}):

\begin{theorem} The polarity map
$$\chi_P \mapsto \chi_{P^\circ}$$
extends to an involution $\Pi:\cC \to \cC$.
\end{theorem}

Now identify $\R^d$ with the open upper hemisphere in $S^d$ by
stereographic projection.  Let $\cD \subset \cC$ be the set of those functions
$f$ such that:
\begin{enumerate}
\item[\bf 1.] $f$ is supported on $\R^d$,
\item[\bf 2.] $f \circ \sigma$ is defined using non-codegenerate
 polytopes, and
\item[\bf 3.] $f \circ \sigma$ is radially left-continuous
as a function on $\R^d$, meaning that for all $v \in \R^d$,
$$\lim_{t \to 1_-} (f \circ \sigma)(tv) = (f \circ \sigma)(v).$$
\end{enumerate}

It is not hard to show that the class $\cD$ is spanned by (the characteristic
functions of) closed convex polytopes with the origin in their interiors.  Thus
$\cD$ is invariant under the polarity involution $\Pi$.

It is also not hard to show that every measure $\mu \in \cA$, if interpreted as
an element of $L^1(\R^d)$, is represented by a unique radially left-continuous
function $f \in \cD$.  This identifies $\cA$ with $\cD$. Again, because $\cD$
is spanned by closed convex polytopes with the origin in their interiors.
Because both maps $\Phi$ and $\Pi$ are the polarity transformation on this
class, the two maps are identified as well.

Note that the class $\cD$ extends to a slightly larger class $\bcD$ spanned by
all convex polytopes $P$ in the closed upper hemisphere in $S^d$ which contain
the origin (not necessarily in the interior).  The class $\bcD$ is also
invariant under the polarity involution $\Pi$.  Indeed $\bcD$ is the closure of
$\cD$ with respect to a natural topology on $\cC$, namely the one induced by
the Hausdorff topology on closed subsets of $S^d$.  Moreover $\Pi$ is
continuous with respect to this topology.  Thus the restriction of $\Pi$ to
$\bcD$ expresses all codegenerate limiting cases of Filliman duality, such as
Lawrence's algorithm.

\acknowledgments

We would like to thank Alexander Barvinok, Jesus De Loera, Yael Karshon, Colin
Rourke, and G\"unter Ziegler for useful discussions.

% \bibliography{mg,gt}

\providecommand{\bysame}{\leavevmode\hbox to3em{\hrulefill}\thinspace}

\end{document}